\documentclass[a4paper,12pt,twoside,leqno,final]{amsart}
\usepackage{amsmath}
\usepackage{amssymb}

\newtheorem{thm}{Theorem}[section]

\newtheorem{cor}[thm]{Corollary}

\theoremstyle{definition}

\newcommand{\C}{{\mathbb C}}
\newcommand{\D}{{\mathbb D}}

\newcommand{\T}{{\mathbb T}}

\newcommand{\bmo}{{\rm BMO}}
\newcommand{\bmoa}{{\rm BMOA}}
\newcommand{\bmona}{{\rm BMO}_{\rm na}}
\newcommand{\vmo}{{\rm VMO}}
\newcommand{\vmoa}{{\rm VMOA}}

\newcommand{\eps}{\varepsilon}

\newcommand{\f}{\frac}
\newcommand{\ov}{\overline}
\newcommand{\al}{\alpha}
\newcommand{\be}{\beta}

\newcommand{\ze}{\zeta}

\newcommand{\si}{\sigma}
\newcommand{\ph}{\varphi}
\newcommand{\om}{\omega}

\newcommand{\const}{\text{\rm const}}

\numberwithin{equation}{section}

\title[Extremal problems in BMO and VMO]
{Extremal problems in BMO and VMO\\ 
involving the Garsia norm}
\author{Konstantin M. Dyakonov}
\address{Departament de Matem\`atiques i Inform\`atica, IMUB, BGSMath, Universitat de Barcelona, Gran Via 585, E-08007 Barcelona, Spain}
\address{ICREA, Pg. Llu\'is Companys 23, E-08010 Barcelona, Spain}
\email{konstantin.dyakonov@icrea.cat}
\thanks{Supported in part by grant PID2021-123405NB-I00 from El Ministerio de Ciencia e Innovaci\'on (Spain) and grant 2021 SGR 00087 from AGAUR (Generalitat de Catalunya).}

\begin{document}
\begin{abstract}
Given an $L^2$ function $f$ on the unit circle $\mathbb T$, we put 
$$\Phi_f(z):=\mathcal P(|f|^2)(z)-|\mathcal Pf(z)|^2,\qquad z\in\mathbb D,$$
where $\mathbb D$ is the open unit disk and $\mathcal P$ is the Poisson integral operator. The Garsia norm $\|f\|_G$ is then defined as $\sup_{z\in\mathbb D}\Phi_f(z)^{1/2}$, and the space ${\rm BMO}$ is formed by the functions $f\in L^2$ with $\|f\|_G<\infty$. If $\|f\|^2_G=\Phi_f(z_0)$ for some point $z_0\in\mathbb D$, then $f$ is said to be a norm-attaining ${\rm BMO}$ function, written as $f\in{\rm BMO}_{\rm na}$. Note that ${\rm BMO}_{\rm na}$ contains ${\rm VMO}$, the space of functions with vanishing mean oscillation. We study, first, the functions $f$ in $L^\infty$ (as well as in $L^\infty\cap{\rm BMO}_{\rm na}$) with the property that $\|f\|_G=\|f\|_\infty$. The analytic case, where $L^\infty$ gets replaced by $H^\infty$, is discussed in more detail. Secondly, we prove that every function $f\in{\rm BMO}_{\rm na}$ with $\|f\|_G=1$ is an extreme point of ${\rm ball}\,({\rm BMO})$, the unit ball of ${\rm BMO}$ with respect to the Garsia norm. This implies that the extreme points of ${\rm ball}\,({\rm VMO})$ are precisely the unit-norm ${\rm VMO}$ functions. As another consequence, we arrive at an amusing \lq\lq geometric" characterization of inner functions.
\end{abstract}

\maketitle

\section{Introduction}

We write $\D$ for the disk $\{z\in\C:|z|<1\}$ and $\T$ for its boundary, the unit circle. The normalized arc length measure on $\T$ is denoted by $m$, and the Lebesgue spaces $L^p:=L^p(\T,m)$ with $0<p\le\infty$ are introduced as usual; the functions that populate them (and those in the other spaces below) are assumed to be complex-valued. With each point $z\in\D$ we associate the {\it harmonic measure} $\om_z$, given by
$$\om_z(E):=\int_E\f{1-|z|^2}{|\ze-z|^2}\,dm(\ze)$$
whenever $E$ is a Lebesgue measurable subset of $\T$. For a function $f\in L^1$, its {\it Poisson integral} $\mathcal Pf$ is defined by
$$\mathcal Pf(z):=\int_\T f(\ze)\,d\om_z(\ze),\qquad z\in\D,$$
so that $\mathcal Pf$ is the harmonic extension of $f$ into $\D$. Finally, we recall that the {\it Hardy space} $H^p$, with $1\le p\le\infty$, is formed by the functions $f\in L^p$ for which $\mathcal Pf$ is analytic on $\D$. Equivalently (see, e.g., \cite[Chapter II]{G}), a function $f$ from $L^p$ belongs to $H^p$ if and only if its {\it Fourier coefficients} 
$$\widehat f(k):=\int_\T\ov\ze^kf(\ze)\,dm(\ze)$$
vanish for all negative indices $k$, that is, 
$$\widehat f(k)=0\quad\text{\rm for }\,\,k=-1,-2,\,\dots$$
When $f\in H^p$, we identify the original function (living a.e. on $\T$) with its analytic extension into $\D$ provided by the Poisson integral, and we denote the latter by $f$ again, rather than by $\mathcal Pf$.

\par Now suppose that $f\in L^2$. We may then consider the nonnegative function $\Phi_f$ on $\D$, built from $f$ by the rule 
\begin{equation}\label{eqn:defcapphi}
\Phi_f(z):=\mathcal P(|f|^2)(z)-|\mathcal Pf(z)|^2,\qquad z\in\D.
\end{equation}
A simple calculation reveals that
$$\Phi_f(z)=\int_\T |f(\ze)-\mathcal Pf(z)|^2\,d\om_z(\ze),\qquad z\in\D,$$
and we can further rewrite the right-hand side as $\mathcal P|f-\mathcal Pf(z)|^2(z)$. Thus, $\Phi_f(z)$ measures the \lq\lq mean deviation from the mean" for $f$ at $z$, where the averages involved are understood as Poisson integrals. We also remark that the function $\Phi_f$ is continuous and superharmonic on $\D$, a fact which follows readily from \eqref{eqn:defcapphi}.

\par We can now define $\bmo=\bmo(\T)$, the space of functions of {\it bounded mean oscillation}, as the set of all functions $f\in L^2$ for which $\Phi_f$ is bounded on $\D$. Its analytic subspace 
$$\bmoa:=\bmo\cap H^2$$
will also be of interest to us. The quantity
$$\|f\|_G:=\sup\left\{\Phi_f(z)^{1/2}:\,z\in\D\right\}$$
(or, strictly speaking, the functional $f\mapsto\|f\|_G$) is known as the {\it Garsia norm} on $\bmo$. Because the functions annihilated by $\|\cdot\|_G$ are the constants, the term \lq\lq seminorm" would actually be more accurate. To make $\|\cdot\|_G$ a true norm, we might therefore agree to identify $\bmo$ functions that differ by a constant. We shall, in fact, adopt this point of view later in the paper (upon special notice), but until then we prefer to treat elements of $\bmo$ as individual functions rather than equivalence classes.

\par We refer to \cite[Chapter VI]{G} or \cite[Chapter X]{K} for various alternative characterizations of (and equivalent norms on) $\bmo$. In addition, these sources offer a systematic presentation of the underlying theory, including the duality between $H^1$ and $\bmoa$, and a great deal more.

\par It is well known that $\bmo$ is contained in every $L^p$ with $0<p<\infty$, while $L^\infty\subset\bmo$. To verify this last inclusion, it suffices to note that 
$$\Phi_f(z)\le\mathcal P(|f|^2)(z)\le\|f\|^2_\infty,\qquad z\in\D,$$
whence
\begin{equation}\label{eqn:inequnifbmo}
\|f\|_G\le\|f\|_\infty
\end{equation}
for any $f\in L^\infty$. (Here, $\|\cdot\|_\infty$ is the usual essential supremum norm on $L^\infty$.) We shall be concerned with cases of equality in \eqref{eqn:inequnifbmo}, and the  functions for which this happens will be referred to as $G$-extremal ones. Thus, a function $f\in L^\infty$ is said to be {\it $G$-extremal} if 
\begin{equation}\label{eqn:equalunifbmo}
\|f\|_G=\|f\|_\infty.
\end{equation}
Roughly speaking, this property means that the function oscillates as much as possible, once its $L^\infty$ norm is fixed. 

\par Among the constant functions, the only $G$-extremal one is the zero function, $f\equiv0$. Less trivial examples can be found among the functions of constant modulus. Consider, for instance, a {\it unimodular} function $f$ on $\T$ (that is, a function $f\in L^\infty$ with $|f|=1$ a.e.) for which $\inf_{z\in\D}|\mathcal Pf(z)|=0$. In this case, \eqref{eqn:defcapphi} reduces to 
$$\Phi_f(z)=1-|\mathcal Pf(z)|^2,\qquad z\in\D,$$
and we infer that $\|f\|_G=1$. Consequently, \eqref{eqn:equalunifbmo} holds true, and $f$ is $G$-extremal. 

\par In particular, we see that every nonconstant inner function is $G$-extremal. It should be recalled that a function $I\in H^\infty$ is said to be {\it inner} if it is unimodular, and that any such $I$ admits a factorization 
$$I=cBS,$$
where $c$ is a constant of modulus $1$, $B$ is a {\it Blaschke product}, and $S$ is a {\it singular 
inner function}; see \cite[Chapter II]{G}. To be more explicit, the factors involved are of the form 
$$B(z)=B_{\{z_j\}}(z):=\prod_j\f{|z_j|}{z_j}\f{z_j-z}{1-\bar z_jz},\qquad z\in\D,$$
where $\{z_j\}\subset\D$ is a sequence---possibly finite or empty---with $\sum_j(1-|z_j|)<\infty$ 
(if $z_j=0$, we put $|z_j|/z_j=-1$), and 
$$S(z)=S_\mu(z):=\exp\left\{-\int_\T\f{\ze+z}{\ze-z}\,d\mu(\ze)\right\},\qquad z\in\D,$$
where $\mu$ is a (nonnegative) singular measure on $\T$. Now, for an inner function $I$, we have 
\begin{equation}\label{eqn:phiinner}
\Phi_I(z)=1-|I(z)|^2,\qquad z\in\D;
\end{equation}
and if $I\ne\const$, then it is also true that $\inf_{z\in\D}|I(z)|=0$, whence $\|I\|_G=1$.

\par For a general analytic function, the presence of an inner factor tends to make the function more oscillatory and thereby increase its Garsia norm. Precisely speaking, for any $F\in H^2$ and any inner function $I$, we have the identity
\begin{equation}\label{eqn:phioutinn}
\Phi_{IF}(z)=\Phi_F(z)+|F(z)|^2(1-|I(z)|^2),\qquad z\in\D.
\end{equation}
One consequence of \eqref{eqn:phioutinn} is the well-known fact that the space $\bmoa$ has the {\it division property} (also called the {\it $f$-property}, see \cite{Hav}): if $F$ and $I$ are as above, and if $IF\in\bmoa$, then $F\in\bmoa$ and $\|F\|_G\le\|IF\|_G$. On the other hand, \eqref{eqn:phioutinn} shows that the set of $G$-extremal functions in $H^\infty$ enjoys the following {\it multiplication property}: if $F(\in H^\infty)$ is $G$-extremal and $I$ is inner, then $IF$ is $G$-extremal. Finally, since for $F\in H^\infty$ and $I$ inner, \eqref{eqn:phioutinn} yields
$$\Phi_{IF}(z)\ge|F(z)|^2(1-|I(z)|^2),\qquad z\in\D,$$
we infer that the condition 
$$\sup\left\{|F(z)|^2(1-|I(z)|^2):\,z\in\D\right\}=\|F\|^2_\infty$$
is sufficient for $IF$ to be $G$-extremal.

\par We go on to recall that $\vmo=\vmo(\T)$, the space of functions of {\it vanishing mean oscillation}, is formed by those $f\in\bmo$ for which 
$$\Phi_f(z)\to0\quad\text{\rm as}\quad|z|\to1^-.$$
The $G$-extremal functions lying in $\vmo$, as well as in 
$$\vmoa:=\vmo\cap H^2,$$ 
will also interest us here. Such functions being bounded, we may equivalently speak of their membership in the space
$$QC:=\vmo\cap L^\infty,$$
whose elements are known as {\it quasicontinuous functions}, and in the corresponding analytic subspace 
$$QA:=\vmoa\cap H^\infty.$$
For basic facts about $QC$ and $QA$, the reader may consult \cite{S} or \cite[Chapter IX]{G}. Among other things, it is shown there that $QC$ (and hence $QA$) is a closed subalgebra of $L^\infty$, and that the inner functions in $QA$ are precisely the rational ones (i.e., finite Blaschke products). 

\par Now, the latter fact provides us with obvious examples of $G$-extremal functions in $QA$. These are of the form $cB$, where $c\in\C$ and $B$ is a nonconstant finite Blaschke product. More interestingly, we shall soon see that this list is exhaustive, meaning that no other examples exist. By contrast, no such rigidity is exhibited by $G$-extremal functions $f$ in $H^\infty$ (let alone $L^\infty$), where it may well happen that $|f|\ne\const$, and moreover, $|f|$ need not obey any additional restrictions whatsoever.

\par It turns out that the rigidity phenomenon just mentioned, which forces $G$-extremal functions to have constant moduli, already manifests itself if the $\vmo$ condition above is replaced by a substantially weaker assumption. Namely, suppose that for a given (possibly unbounded) function $f\in\bmo$ we can find a point $z_0\in\D$ such that 
$$\sup\left\{\Phi_f(z):\,z\in\D\right\}=\Phi_f(z_0)$$
(or, equivalently, $\|f\|^2_G=\Phi_f(z_0)$). If so, we say that $f$ is a {\it norm-attaining} $\bmo$ function, and we denote the set of such functions by $\bmona$. This class will turn up repeatedly in what follows.

\par It should be noted that 
\begin{equation}\label{eqn:bmonaincl}
\vmo\subset\bmona\subset\bmo,
\end{equation}
the three sets being all different. In particular, to check the left-hand inclusion in \eqref{eqn:bmonaincl}, we fix a function $f\in\vmo$ and observe that $\Phi_f$ extends continuously to $\D\cup\T$, once we define it to be zero on $\T$, so $\Phi_f$ must attain its maximum value at some point of $\D$. Furthermore, given an inner function $I$, identity \eqref{eqn:phiinner} tells us that $I\in\bmona$ if and only if either $I=\const$ or $I$ has a zero in $\D$. Recalling that every inner function is in $\bmo$, while no non-rational inner function is in $\vmo$, we see that both inclusions in \eqref{eqn:bmonaincl} are indeed proper. To get an example of a non-analytic function $f\in\bmona\setminus\vmo$, one can take $f=\chi_E$, the characteristic function of a measurable set $E\subset\T$ with $0<m(E)<1$. In this case, 
$$\Phi_f(z)=\om_z(E)\cdot\left(1-\om_z(E)\right),\qquad z\in\D,$$
and the required properties are easy to verify.

\par We also remark that the set $\bmona$ is nonlinear, i.e., it does not form a vector space. Indeed, the singular inner function $S(z):=\exp\{(z+1)/(z-1)\}$ can be written as 
$$S=(1-z)S+zS$$
and both terms on the right-hand side are in $\bmona$, whereas $S$ is not.

\par Our results are stated in Section 2 below. These chiefly concern the moduli of $G$-extremal functions on $\T$, the canonical factorization of such functions (in the analytic setting), and the geometry of the unit balls in the spaces $\bmo$ and $\vmo$, always endowed with the Garsia norm $\|\cdot\|_G$. We also mention a new \lq\lq geometric" characterization of inner functions that comes out in this context. The proofs are given in Sections 3 and 4. Finally, Section 5 (which can be viewed as a kind of addendum) contains an explicit---but somewhat technical---construction of a $G$-extremal outer function in $H^\infty$.

\par After this work had been completed, Kristian Seip kindly drew my attention to the recent paper \cite{BS}, where a related extremal problem was treated in connection with Hankel operators. There, of concern were the functions $f\in H^\infty$ that satisfy
$$\|f\|_\infty=\inf\left\{\|f-\ov z\ov h\|_\infty:\,h\in H^\infty\right\},$$
a property reminiscent of \eqref{eqn:equalunifbmo}. The class that arises turns out to have quite a bit in common with that of $G$-extremal functions lying in $H^\infty$. In particular, the similarity becomes apparent upon juxtaposing Theorem 3 from \cite{BS} with our Corollary \ref{cor:factdiskalg}, to be found in the next section.

\section{Statement of results}

We begin by looking at the moduli of $G$-extremal functions. Among other things, we show (see part (a) of the theorem below) that every nonnegative $L^\infty$ function on $\T$ is expressible as $|f|$ for some $G$-extremal function $f$ from a fairly small subalgebra of $L^\infty$, not just from $L^\infty$ itself. To be precise, the subalgebra in question is $H^\infty+C$, the linear hull of $H^\infty$ and $C:=C(\T)$. It is well known that $H^\infty+C$ is actually the smallest closed subalgebra of $L^\infty$ that contains $H^\infty$ properly, while $QC$ is the largest self-adjoint subalgebra of $H^\infty+C$; see \cite[Chapter IX]{G} or \cite{S}. 

\begin{thm}\label{thm:modgextr} Let $\eta$ be a nonnegative function in $L^\infty$, $\eta\not\equiv0$.
	
\par{\rm (a)} We have then 
\begin{equation}\label{eqn:phimodf}
\eta=|f|\quad\text{a.e. on }\T
\end{equation}	
for some $G$-extremal function $f\in H^\infty+C$. 
\par{\rm (b)} In order that $\eta$ be writable in the form \eqref{eqn:phimodf} with some $G$-extremal function $f\in H^\infty$, it is necessary and sufficient that
\begin{equation}\label{eqn:logintconv}
\int_\T\log\eta\,dm>-\infty.
\end{equation}
\par{\rm (c)} In order that $\eta$ be writable in the form \eqref{eqn:phimodf} with some $G$-extremal function $f\in\bmona\cap L^\infty$, it is necessary and sufficient that $\eta\equiv\const$. 
\end{thm}

\par Of course, the necessity in (b) is trivial, since \eqref{eqn:logintconv} characterizes the moduli of non-null $H^\infty$ functions (see, e.g. \cite[Chapter II]{G}). Also obvious is the sufficiency part in (c). Indeed, given a constant $c>0$, the function $f(z)=cz$ has modulus $c$ on $\T$, is $G$-extremal and belongs to $\bmona\cap L^\infty$. This observation actually shows that statement (c) remains valid whenever the class $\bmona\cap L^\infty$ is replaced by a subset thereof that contains the identity function $z$ (among such subsets we mention $QC$, $QA$ and $\bmona\cap H^\infty$).

\par Combining part (c) of Theorem \ref{thm:modgextr} with our preliminary observations on inner functions from Section 1, we arrive at the following result. 

\begin{cor}\label{cor:ginn} {\rm (1)} The general form of a $G$-extremal function in $\bmona\cap H^\infty$ is $cI$, where $c\in\C$ and $I$ is an inner function that has a zero in $\D$.
\par{\rm (2)} The general form of a $G$-extremal function in $QA$ is $cB$, where $c\in\C$ and $B$ is a nonconstant finite Blaschke product. 
\end{cor}

\par The deeper part of Theorem \ref{thm:modgextr}, namely statement (a) therein, hinges on a theorem of Axler (see \cite{Ax77}) saying that every $L^\infty$ function is writable as $f/B$, where $f\in H^\infty+C$ and $B$ is a Blaschke product. We actually need a refinement of Axler's result, stated as Theorem \ref{thm:imprax} below, where $f$ is additionally taken to be $G$-extremal.

\begin{thm}\label{thm:imprax} Let $\ph\in L^\infty$. Then there exist a Blaschke product $B$ and a $G$-extremal function $f\in H^\infty+C$ such that $B\ph=f$.
\end{thm}

\medskip\noindent{\it Remark.} If we assume that $\ph\in H^\infty$, then the function $f$ provided by Theorem \ref{thm:imprax} will also be in $H^\infty$, not just in $H^\infty+C$; indeed, the product $B\ph$ is then obviously in $H^\infty$. This \lq\lq analytic version" of the above theorem is, however, essentially simpler and can be proved without recourse to Axler's original result (while our proof of Theorem \ref{thm:imprax} does rely on it). In Section 3, we use that alternative---and more direct---approach to establish part (b) of Theorem \ref{thm:modgextr}.

\medskip Before proceeding further, we recall that an analytic function $F$ on $\D$ is said to be {\it outer} if $F$ is zero-free and the harmonic function $\log|F(z)|$ is the Poisson integral of an $L^1$ function on $\T$. Up to a constant factor of modulus $1$, such a function has the form 
\begin{equation}\label{eqn:outfuneta}
\mathcal O_\eta(z):=\exp\left\{\int_\T\f{\ze+z}{\ze-z}\log\eta(\ze)\,dm(\ze)\right\},\qquad z\in\D,
\end{equation}
where $\eta$ is some nonnegative function on $\T$ with $\log\eta\in L^1$. Since 
$$\left|\mathcal O_\eta(z)\right|=\exp\left\{\mathcal P(\log\eta)(z)\right\},\qquad z\in\D,$$
the boundary values of $\mathcal O_\eta$ (understood as nontangential limits) satisfy $|\mathcal O_\eta|=\eta$ a.e. on $\T$. 

\par It is a classical result (see, e.g., \cite[Chapter II]{G}) that the moduli of non-null $H^p$ functions, and hence the outer functions $\mathcal O_\eta$ in $H^p$, are characterized by the conditions $\eta\in L^p$ and $\log\eta\in L^1$; equivalently, the former condition should be coupled with \eqref{eqn:logintconv}. The $\bmo$ situation is subtler, and we refer to \cite{DRMI} for a description of the outer functions $\mathcal O_\eta$ lying in $\bmoa$. In \cite{DIUMJ}, the membership of $\mathcal O_\eta$ in $QA$ was discussed (as part of a more general problem) and the appropriate criterion was furnished.

\par Going back to Theorem \ref{thm:modgextr}, part (b), we may rephrase that statement by saying that every outer function $F\in H^\infty$ occurs in the canonical (inner-outer) factorization of some $G$-extremal $H^\infty$ function. In other words, any such $F$ can be multiplied by a suitable inner factor $I$ so that the product $IF$ becomes $G$-extremal. It seems natural to ask which outer functions (if any) are $G$-extremal in their own right, no help from the inner factor being needed. Such outer functions do exist, but they are necessarily pretty discontinuous (not even in $QA$), and we characterize them as follows.

\begin{thm}\label{thm:charoutergextr} Let $\eta\in L^\infty$ be a nonnegative function with $\int_\T\log\eta\,dm>-\infty$, and let $\mathcal O_\eta$ be the outer function given by \eqref{eqn:outfuneta}. Then $\mathcal O_\eta$ is $G$-extremal if and only if there exists a sequence $\{z_n\}\subset\D$ such that $\mathcal P\eta(z_n)\to\|\eta\|_\infty$ and 
$$\mathcal P\left(\log\f1\eta\right)(z_n)\to\infty.$$
\end{thm}

\par This characterization is fairly easy to establish, but it may not be immediately clear how to produce a function $\eta$ that meets the criterion. An explicit example is constructed in Section 5 below. 

\par Now consider an outer function $F\in QA$. It follows from Corollary \ref{cor:ginn}, part (2), that $F$ is never $G$-extremal. At the same time, we already know that the $G$-extremality can be achieved by passing from $F$ to $IF$ if the inner factor $I$ is chosen appropriately. Our next result characterizes the inner functions $I$ that are eligible in this sense. The criterion that arises, as given in Theorem \ref{thm:factgextr} below, might be contrasted with that in \cite[Theorem 4]{DAJM}, where the canonical factorization of a $QA$ function was treated. There, the outer factor was required to compensate for the oscillatory behavior of the inner one, while in the current setting the inner factor is needed to make the product oscillatory enough.

\begin{thm}\label{thm:factgextr} Suppose that $F$ is an outer function in $QA$ and $I$ is an inner function. Then $IF$ is $G$-extremal if and only if there exists a sequence $\{z_n\}\subset\D$ such that $|F(z_n)|\to\|F\|_\infty$ and $I(z_n)\to0$.
\end{thm}

\par This criterion takes an even simpler form if $F$ is assumed to be in the {\it disk algebra} $\mathcal A:=H^\infty\cap C(\T)$, not just in $QA$. To state this new version, we first introduce a bit of notation. For $F\in\mathcal A$, we define 
$$\mathcal M(F):=\left\{z\in\D\cup\T:\,|F(z)|=\|F\|_\infty\right\}.$$
Clearly, this set either coincides with $\D\cup\T$ (which happens precisely when $F$ is constant) or is contained in $\T$. Also, for an inner function $I$, we consider the set
$$\si(I):=\left\{z\in\D\cup\T:\,\liminf_{\D\ni w\to z}|I(w)|=0\right\},$$
known as the {\it spectrum} of $I$. Equivalently, we have 
$$\si(I)=\overline{\mathcal Z_I\cup\mathcal S(\mu_I)},$$
where $\mathcal Z_I$ is the zero set of $I$ in $\D$, $\mathcal S(\mu_I)$ is the support of the associated singular measure $\mu=\mu_I$ on $\T$, and the bar stands for closure. 

\begin{cor}\label{cor:factdiskalg} Given an outer function $F\in\mathcal A$ and an inner function $I$, the product $IF$ is $G$-extremal if and only if $\mathcal M(F)\cap\si(I)\ne\emptyset$.
\end{cor}

\par To deduce this from Theorem \ref{thm:factgextr}, suppose first that $\mathcal M(F)\cap\si(I)$ contains a point $z^*$. Then there is a sequence $\{z_n\}\subset\D$ such that $z_n\to z^*$ and $I(z_n)\to0$. We also have then $|F(z_n)|\to|F(z^*)|=\|F\|_\infty$, so $\{z_n\}$ has the properties specified in Theorem \ref{thm:factgextr}. Conversely, if $\{z_n\}\subset\D$ is a sequence coming from Theorem \ref{thm:factgextr}, then a suitable subsequence $\{z_{n_k}\}$ is convergent (to a point in $\D\cup\T$) and the limit point is actually in $\mathcal M(F)\cap\si(I)$.

\medskip\noindent{\it Remark.} The above reasoning, and hence Corollary \ref{cor:factdiskalg}, remains valid if we merely assume that the outer function $F$ is in $H^\infty$ and $|F|$ is continuous on $\D\cup\T$. (These conditions are known to imply that $F\in QA$; see \cite[Theorem 3]{Sar}.) The current version of the corollary is, however, good enough to produce numerous examples of analytic $G$-extremal functions with \lq\lq nice" outer factors.

\medskip Yet another---not totally unrelated---problem to be treated here concerns the extreme points of the unit ball in $\bmo$ (and in $\vmo$). To keep on the safe side, we emphasize that we no longer restrict ourselves to bounded functions. Given a Banach space $X$, we write ${\rm ball}\,(X)$ for its closed unit ball, so that 
$$\text{\rm ball}\,(X):=\{x\in X:\,\|x\|\le1\}.$$
As usual, we say that an element $x$ of $\text{\rm ball}\,(X)$ is an {\it extreme point} thereof if it is not expressible in the form $x=\f12(u+v)$ with two distinct points $u,v\in\text{\rm ball}\,(X)$. Of course, every extreme point of $\text{\rm ball}\,(X)$ has norm $1$. In what follows, we are chiefly interested in the case where $X$ is either $\bmo$ or $\vmo$. Throughout, both spaces are endowed with the Garsia norm $\|\cdot\|_G$ and their unit balls are defined accordingly. Also, to make sure that $\|\cdot\|_G$ is actually a norm---rather than a seminorm---on $\bmo$ and/or $\vmo$, we shall henceforth regard each of these as a quotient space modulo constants. Thus, an element $f$ of $\bmo$ or $\vmo$, even if referred to as a function, should be viewed as a coset (i.e., equivalence class) of the form $f+\C$. 

\par For other types of norms on $\bmo$ and $\vmo$, the problem of characterizing the corresponding extreme points was considered by Axler and Shields; see \cite{AS}. (The norms involved were defined in terms of averages over arcs.) The common feature between their work and ours is that the $\vmo$ case admits a neat solution, while for $\bmo$ only partial results are provided. 

\par First, we come up with the following sufficient condition.

\begin{thm}\label{thm:extpoibmo} If $f\in\bmona$ and $\|f\|_G=1$, then $f$ is an extreme point of ${\rm ball}\,(\bmo)$.
\end{thm}

\par Clearly, a similar result holds in the analytic setting where $\bmo$ and $\bmona$ get replaced by $\bmoa$ and ${\rm BMOA}_{\rm na}:=\bmona\cap H^2$, respectively. 
\par On the other hand, since $\vmo\subset\bmona$, it follows from Theorem \ref{thm:extpoibmo} that every unit-norm $\vmo$ function is an extreme point of ${\rm ball}\,(\bmo)$ and hence also of ${\rm ball}\,(\vmo)$, a fact that makes the geometry of the latter set remarkably simple. We summarize this observation as a corollary below.

\begin{cor}\label{cor:extpoivmo} The extreme points of ${\rm ball}\,(\vmo)$ are precisely the functions $f\in\vmo$ with $\|f\|_G=1$.
\end{cor}

\par We mention in passing that Theorem \ref{thm:extpoibmo} and Corollary \ref{cor:extpoivmo} have natural analogs for some other classical spaces that admit a Garsia-type norm. Consider, by way of example, the analytic {\it Lipschitz space} $A^\al$ with $0<\al<\f12$. By definition, $A^\al$ is the set of all functions $f\in H^\infty$ that satisfy 
\begin{equation}\label{eqn:lipcondalpha}
|f(z)-f(w)|\le C_f|z-w|^\al\qquad(z,w\in\D)
\end{equation}
with some constant $C_f\ge0$. While $A^\al$ is usually normed by letting $\|f\|_{A^\al}$ be the best value of $C_f$ in \eqref{eqn:lipcondalpha}, it is known (see \cite[Lemma 1]{DMMJ} or \cite[Theorem 1]{DActa}) that an equivalent norm is provided by the Garsia-type quantity
$$\|f\|_{G,\al}:=\sup\left\{\f{\sqrt{\Phi_f(z)}}{(1-|z|)^\al}:\,z\in\D\right\},$$
as long as $0<\al<\f12$. (To check that the threshold exponent $\f12$ is sharp, take $f$ to be the identity function $z$.) Furthermore, if we define the analytic {\it little Lipschitz space} $A_0^\al$, with $\al$ as above, by means of the appropriate \lq\lq little oh" version of \eqref{eqn:lipcondalpha}, then we have 
$$A_0^\al=\{f\in A^\al:\,\Phi_f(z)=o((1-|z|)^{2\al})\,\,\text{\rm as }|z|\to1^-\}.$$
Both spaces, $A^\al$ and $A_0^\al$, are endowed below with the Garsia-type norm $\|\cdot\|_{G,\al}$, their unit balls being defined accordingly, and both are viewed as quotient spaces modulo constants. 

\par Finally, we introduce the set $A_{\rm na}^\al$ of norm-attaining $A^\al$ functions (with respect to $\|\cdot\|_{G,\al}$). Given a function $f$ in $A^\al$, with $0<\al<\f12$, we write $f\in A_{\rm na}^\al$ if there is a point $z_0\in\D$ such that 
$$\|f\|_{G,\al}=\f{\sqrt{\Phi_f(z_0)}}{(1-|z_0|)^\al}.$$
It should be noted that $A_0^\al\subset A_{\rm na}^\al$. The counterparts of Theorem \ref{thm:extpoibmo} and Corollary \ref{cor:extpoivmo} that arise in this setting read as follows. 

\begin{thm}\label{thm:extpoilip} Let $0<\al<\f12$. If $f\in A_{\rm na}^\al$ and $\|f\|_{G,\al}=1$, then $f$ is an extreme point of ${\rm ball}\,(A^\al)$.
\end{thm}

\begin{cor}\label{cor:extpoilitlip} For $0<\al<\f12$, the extreme points of ${\rm ball}\,(A_0^\al)$ are precisely the functions $f\in A_0^\al$ with $\|f\|_{G,\al}=1$.
\end{cor}

\par Going back to the extreme points of the unit ball in $\bmo$, or in $\bmoa$, we do not know whether the sufficient condition coming from Theorem \ref{thm:extpoibmo} is actually necessary. (A similar question can be asked in connection with Theorem \ref{thm:extpoilip}.) Anyhow, in contrast to the $\vmo$ situation, the boundary sphere has this time a rich supply of non-extreme points, as the following result shows. 

\begin{thm}\label{thm:nonextpoibmo} {\rm (a)} Given a non-unimodular function $\ph\in L^\infty$ with $\|\ph\|_\infty=1$, there exists a Blaschke product $B$ such that $\|B\ph\|_G=1$ and $B\ph$ is a non-extreme point of ${\rm ball}\,(\bmo)$.
\par{\rm (b)} Given a function $\ph\in H^\infty$ with $\|\ph\|_\infty=1$ for which
\begin{equation}\label{eqn:logintegconv}
\int_\T\log(1-|\ph|)\,dm>-\infty,
\end{equation}
there exists a Blaschke product $B$ such that $\|B\ph\|_G=1$ and $B\ph$ is a non-extreme point of ${\rm ball}\,(\bmoa)$.
\par{\rm (c)} Moreover, if $\ph$ is a function in $L^\infty$ (resp., in $H^\infty$) with $\|\ph\|_\infty=1$ which is non-unimodular (resp., which satisfies \eqref{eqn:logintegconv}) and if $B$ is any Blaschke product, then $B\ph$ is a non-extreme point of ${\rm ball}\,(\bmo)$ (resp., of ${\rm ball}\,(\bmoa)$).
\end{thm}

\par It should be recalled that the extreme points of ${\rm ball}\,(L^\infty)$ are precisely the unimodular functions, while those of ${\rm ball}\,(H^\infty)$ are the functions $\ph\in H^\infty$ with $\|\ph\|_\infty=1$ for which \eqref{eqn:logintegconv} fails; see \cite[Section V]{dLR} or \cite[Chapter 9]{Hof}. Thus, our method of constructing a unit-norm non-extreme point of ${\rm ball}\,(\bmo)$ (resp., ${\rm ball}\,(\bmoa)$), as employed in Theorem \ref{thm:nonextpoibmo} above, consists in multiplying a unit-norm non-extreme function $\ph$ from ${\rm ball}\,(L^\infty)$ (resp., ${\rm ball}\,(H^\infty)$) by a suitable Blaschke product $B$; the latter is needed to make the resulting function $G$-extremal.

\par Finally, much in the same spirit, we arrive at an amusing characterization of the inner functions among all unit-norm elements of $H^\infty$. 

\begin{thm}\label{thm:charint} Suppose that $h\in H^\infty$ and $\|h\|_\infty=1$. The following conditions are equivalent: 
\par{\rm (i)} $h$ is an inner function. 
\par{\rm (ii)} For each nonconstant Blaschke product $B$, the function $Bh$ is an extreme point of ${\rm ball}\,(\bmo)$. 
\par{\rm (iii)} There exists a Blaschke product $B$ such that $Bh$ is an extreme point of ${\rm ball}\,(\bmo)$.
\end{thm}

\par Note that, even though the functions involved (i.e., $h$ and $B$) are in $H^\infty$, it is the \lq\lq non-analytic" space $\bmo$, rather than $\bmoa$, that figures in conditions (ii) and (iii) above. On the other hand, this last theorem is vaguely reminiscent of another \lq\lq geometric" result, due to Cima and Thomson, which identifies the inner functions with the so-called strongly extreme points of ${\rm ball}\,(H^\infty)$. A proof of this fact---and the definition of a strongly extreme point---can be found in \cite{CT}; see also \cite{DAFM} and \cite{MvR} for generalizations.

\par We now turn to the proofs of our results. These are provided in the next two sections.

\section{Proofs of Theorems \ref{thm:modgextr}, \ref{thm:imprax}, \ref{thm:charoutergextr} and \ref{thm:factgextr}}

We prove Theorem \ref{thm:imprax} first. Once this is done, the bulk of Theorem \ref{thm:modgextr} will be established as well.

\medskip\noindent{\it Proof of Theorem \ref{thm:imprax}.} We may discard the trivial case where $\ph\equiv0$ and assume, without any loss of generality, that $\|\ph\|_\infty=1$. By a result of Axler (see \cite[Theorem 1]{Ax77}), there exist a Blaschke product $b$ and a function $\psi\in H^\infty+C$ such that 
\begin{equation}\label{eqn:bfipsi}
b\ph=\psi. 
\end{equation}
Since $|\ph|^2$ is a unit-norm function in $L^\infty$, we can find a sequence $\{z_n\}\subset\D$ with $|z_n|\to1$ for which
\begin{equation}\label{eqn:poitoone}
\mathcal P\left(|\ph|^2\right)(z_n)\to1.
\end{equation}
Passing to a subsequence if necessary, we may assume in addition that $z_n\to\ze_0$ for some $\ze_0\in\T$ and that $\sum_n(1-|z_n|)<\infty$. We then define $\mathcal B$ to be the Blaschke product with zeros $\{z_n\}$. Next, we introduce the Blaschke product $B:=\mathcal Bb$ and the function $f:=\mathcal B\psi$, whereupon \eqref{eqn:bfipsi} takes the form $B\ph=f$. It follows, in particular, that 
\begin{equation}\label{eqn:modequal}
|f|=|\ph|\quad\text{\rm a.e. on }\T.
\end{equation}
\par The proof will be complete once we show that $f\in H^\infty+C$ and $\|f\|_G=1$. (The latter condition will imply that $f$ is $G$-extremal, since $\|f\|_\infty=\|\ph\|_\infty=1$.) Recalling that $\psi\in H^\infty+C$, we can write 
\begin{equation}\label{eqn:decomp}
\psi=g+h\quad\text{\rm with}\quad g\in H^\infty,\,\,h\in C,
\end{equation}
where the continuous function $h$ is taken to satisfy
\begin{equation}\label{eqn:adprop}
h(\ze_0)=0.
\end{equation}
(Such a choice is, of course, possible: given any decomposition \eqref{eqn:decomp}, the additional property \eqref{eqn:adprop} is achieved upon replacing $g$ and $h$ by $g+h(\ze_0)$ and $h-h(\ze_0)$, respectively.) We have then 
\begin{equation}\label{eqn:vbgbh}
f=\mathcal B\psi=\mathcal Bg+\mathcal Bh
\end{equation}
and we observe that $\mathcal Bg\in H^\infty$, while $\mathcal Bh\in C$; this last fact is ensured by \eqref{eqn:adprop}, since $\ze_0$ is the only discontinuity point for $\mathcal B$ on $\T$. Thus, \eqref{eqn:vbgbh} tells us that $f\in H^\infty+C$. 

\par Furthermore, for $z\in\D$, \eqref{eqn:vbgbh} yields
$$\mathcal Pf(z)=\mathcal B(z)g(z)+\mathcal P(\mathcal Bh)(z).$$
Consequently, 
\begin{equation*}
\begin{aligned}
\Phi_f(z)&=\mathcal P\left(|f|^2\right)(z)-|\mathcal Pf(z)|^2\\
&=\mathcal P\left(|\ph|^2\right)(z)-\left|\mathcal B(z)g(z)+\mathcal P(\mathcal Bh)(z)\right|^2,
\end{aligned}
\end{equation*}
where we have also used \eqref{eqn:modequal}. Letting $z=z_n$, we obtain
\begin{equation}\label{eqn:phivzn}
\Phi_f(z_n)=\mathcal P\left(|\ph|^2\right)(z_n)-\left|\mathcal P(\mathcal Bh)(z_n)\right|^2.
\end{equation}
Finally, we recall \eqref{eqn:poitoone} and note that $\mathcal P(\mathcal Bh)(z_n)\to0$ (because $z_n\to\ze_0$, while $\mathcal Bh$ is a continuous function on $\T$ that vanishes at $\ze_0$). In light of these facts, \eqref{eqn:phivzn} shows that $\Phi_f(z_n)\to1$. This in turn implies that $\|f\|_G=1$, as desired. \qed

\bigskip\noindent{\it Proof of Theorem \ref{thm:modgextr}.} Part (a) is immediate from Theorem \ref{thm:imprax}. Indeed, applying that theorem with $\ph=\eta$, we arrive at a $G$-extremal 
function $f\in H^\infty+C$ that satisfies $|f|=\eta$ and therefore does the job.

\par We shall assume throughout the rest of the proof that $\|\eta\|_\infty=1$. Turning to part (b), we only have to prove the sufficiency of \eqref{eqn:logintconv}. First of all, \eqref{eqn:logintconv} allows us to define the outer function $\mathcal O_\eta$ in accordance with \eqref{eqn:outfuneta}, and we put $F:=\mathcal O_\eta$. Thus, in particular, $F\in H^\infty$ and $|F|=\eta$ a.e. on $\T$. Now, one quick way of producing a $G$-extremal function $f\in H^\infty$ with $|f|=\eta$ would be to invoke Theorem \ref{thm:imprax} once again (see also the remark attached to it), this time with $\ph=F$. An alternative route, taken below, leads to a more explicit construction of a Blaschke product $B$ for which the product $BF(=f)$ is $G$-extremal.

\par Since $\eta$, and hence $\eta^2$, has norm $1$ in $L^\infty$, there is a sequence $\{z_n\}\subset\D$ with $|z_n|\to1$ satisfying
\begin{equation}\label{eqn:poiphitoone}
\mathcal P\left(\eta^2\right)(z_n)\to1.
\end{equation}
Assuming in addition that $\sum_n(1-|z_n|)<\infty$ (which is achieved by passing to a subsequence if necessary), we write $B$ for the Blaschke product with zeros $\{z_n\}$; then we put $f:=BF$. It is clear that $f\in H^\infty$, $\|f\|_\infty=1$, and that \eqref{eqn:phimodf} is fulfilled. Besides, for $z\in\D$, we have 
\begin{equation*}
\begin{aligned}
\Phi_f(z)&=\mathcal P\left(|f|^2\right)(z)-|f(z)|^2\\
&=\mathcal P\left(\eta^2\right)(z)-|B(z)|^2|F(z)|^2,
\end{aligned}
\end{equation*}
whence 
$$\Phi_f(z_n)=\mathcal P\left(\eta^2\right)(z_n)$$
for all $n$. Recalling \eqref{eqn:poiphitoone}, we conclude that 
$$\sup_{z\in\D}\Phi_f(z)\left(=\|f\|^2_G\right)=1$$
and so $f$ is $G$-extremal, as required. 

\par As regards part (c), this time it is the necessity we need to verify. Suppose $f$ is a $G$-extremal function in $\bmona\cap L^\infty$ that satisfies \eqref{eqn:phimodf}. Together with the assumption that $\|\eta\|_\infty=1$, this means that there is a point $z_0\in\D$ for which $$\Phi_f(z_0)=1\left(=\|f\|^2_G\right),$$
or equivalently,
\begin{equation}\label{eqn:smtheqtoone}
\mathcal P\left(\eta^2\right)(z_0)-|\mathcal Pf(z_0)|^2=1.
\end{equation}
Because the quantity 
\begin{equation}\label{eqn:skovoroda}
\mathcal P\left(\eta^2\right)(z_0):=\int_\T\eta^2\,d\omega_{z_0}
\end{equation}
is obviously bounded by $1$, identity \eqref{eqn:smtheqtoone} may only hold if 
\begin{equation}\label{eqn:solonka}
\mathcal P\left(\eta^2\right)(z_0)=1,
\end{equation}
while $\mathcal Pf(z_0)=0$. This in turn implies that $\eta=1$ a.e. on $\T$. (To see why, combine \eqref{eqn:skovoroda} and \eqref{eqn:solonka} with the facts that $0\le\eta\le1$ and $\omega_{z_0}(\T)=1$.) We thus arrive at the desired conclusion that $\eta$ is constant. \qed

\bigskip\noindent{\it Proof of Theorem \ref{thm:charoutergextr}.} Let $F:=\mathcal O_\eta$. Assume first that $\|\eta\|_\infty=1$ (and hence $\|F\|_\infty=1$). In order that $F$ be $G$-extremal, it is necessary and sufficient that there exist a sequence $\{z_n\}\subset\D$ such that 
$$\Phi_F(z_n)\to1,$$
or equivalently, 
\begin{equation}\label{eqn:lozhka}
\mathcal P\left(\eta^2\right)(z_n)-\exp\left\{2\mathcal P(\log\eta)(z_n)\right\}\to1.
\end{equation}
Because $\mathcal P\left(\eta^2\right)(z_n)\le1$, condition \eqref{eqn:lozhka} means that 
\begin{equation}\label{eqn:vilka}
\mathcal P\left(\eta^2\right)(z_n)\to1
\end{equation}
and
\begin{equation}\label{eqn:nozhik}
\exp\left\{2\mathcal P(\log\eta)(z_n)\right\}\to0.
\end{equation}
Furthermore, using the fact that $0\le\eta\le1$ and the Cauchy--Schwarz inequality, we get 
$$\mathcal P\left(\eta^2\right)(z_n)\le\mathcal P\eta(z_n)\le
\left\{\mathcal P\left(\eta^2\right)(z_n)\right\}^{1/2},$$
which allows us to rewrite \eqref{eqn:vilka} in the form 
\begin{equation}\label{eqn:vilkabis}
\mathcal P\eta(z_n)\to1.
\end{equation}
On the other hand, we can rephrase \eqref{eqn:nozhik} as 
\begin{equation}\label{eqn:nozhikbis}
\mathcal P\left(\log\f1\eta\right)(z_n)\to\infty.
\end{equation}
To summarize, our criterion of $G$-extremality amounts under the current hypotheses to the existence of a sequence $\{z_n\}\subset\D$ that makes \eqref{eqn:vilkabis} and \eqref{eqn:nozhikbis} true.
\par Now, to remove the assumption that $\|\eta\|_\infty=1$, we only need to apply the above criterion with $\eta/\|\eta\|_\infty$ in place of $\eta$. In doing so, we may however leave \eqref{eqn:nozhikbis} intact; indeed, this last condition remains unchanged upon replacing $\eta$ by $c\eta$ with any constant $c>0$.\qed

\bigskip\noindent{\it Proof of Theorem \ref{thm:factgextr}.} We shall assume, without losing anything in generality, that $\|F\|_\infty=1$. 
\par Suppose, to begin with, that $|F|=1$ a.e. on $\T$. We have then $F\equiv c$ for some constant $c\in\C$ with $|c|=1$. In that case, the product $IF$ is $G$-extremal whenever $I$ is a nonconstant inner function, and the desired equivalence relation follows easily. (Indeed, for $I$ inner, the existence of a sequence $\{z_n\}\subset\D$ with $I(z_n)\to0$ means precisely that $I\ne\const$, while the condition $|F(z_n)|\to1$ is automatic since $|F|=1$ on $\D$.) 
\par Now suppose that $|F|$ is nonconstant on $\T$. Given an inner function $I$, the product $IF=:f$ is $G$-extremal if and only if there is a sequence $\{z_n\}\subset\D$ such that $\Phi_f(z_n)\to1$. Furthermore, we necessarily have $|z_n|\to1$. (Otherwise, $\{z_n\}$ would have a subsequence converging to some point $z^*\in\D$, and it would follow that $\Phi_f(z^*)=1$, so that $f\in\bmona\cap H^\infty$. By Theorem \ref{thm:modgextr}, part (c), this would contradict the assumption that $|f|(=|F|)$ is nonconstant on $\T$.) Using the identity 
$$\Phi_f(z_n)=\Phi_F(z_n)+|F(z_n)|^2(1-|I(z_n)|^2)$$
together with the fact that $\Phi_F(z_n)\to0$ (which holds because $|z_n|\to1$ and $F\in QA$), we now rewrite the condition $\Phi_f(z_n)\to1$ in the form 
\begin{equation}\label{eqn:skalka}
|F(z_n)|^2(1-|I(z_n)|^2)\to1.
\end{equation}
In view of the obvious inequalities 
$$|F(z_n)|^2(1-|I(z_n)|^2)\le|F(z_n)|^2\le1,$$
we can finally rephrase \eqref{eqn:skalka} by saying that $|F(z_n)|\to1$ and $I(z_n)\to0$, as required.
\qed

\section{Proofs of Theorems \ref{thm:extpoibmo}, \ref{thm:extpoilip}, \ref{thm:nonextpoibmo} and \ref{thm:charint}}

\noindent{\it Proof of Theorem \ref{thm:extpoibmo}.} Let $f\in\bmona$ be a function with $\|f\|_G=1$. Thus, for some $z_0\in\D$, we have 
\begin{equation}\label{eqn:phifeqone}
\Phi_f(z_0)=1.
\end{equation}
To prove that $f$ (or rather the coset $f+\C$) is an extreme point of ${\rm ball}\,({\bmo})$, we take an arbitrary function $g\in\bmo$ and show that the inequalities 
\begin{equation}\label{eqn:fpmg}
\|f+g\|_G\le1\quad\text{\rm and}\quad\|f-g\|_G\le1
\end{equation}
may only hold if $g$ is constant (so that the coset $g+\C$ is null). 
\par For a fixed $z\in\D$, consider the functions 
$$\mathcal F_z(\ze):=f(\ze)-\mathcal Pf(z)$$
and
$$\mathcal G_z(\ze):=g(\ze)-\mathcal Pg(z),$$
defined for almost all $\ze\in\T$. Integrating the identity 
$$2|\mathcal F_z|^2+2|\mathcal G_z|^2=|\mathcal F_z+\mathcal G_z|^2+|\mathcal F_z-\mathcal G_z|^2$$
over $\T$ against the harmonic measure $d\om_z$, we obtain 
\begin{equation}\label{eqn:pmphiphi}
2\Phi_f(z)+2\Phi_g(z)=\Phi_{f+g}(z)+\Phi_{f-g}(z).
\end{equation}
Because 
$$\Phi_{f\pm g}(z)\le\|f\pm g\|^2_G\le1$$
as long as \eqref{eqn:fpmg} holds, the right-hand side of \eqref{eqn:pmphiphi} is bounded by $2$, and we infer that 
$$\Phi_g(z)\le1-\Phi_f(z)$$
for all $z\in\D$. In view of \eqref{eqn:phifeqone}, this implies that 
\begin{equation}\label{eqn:phigeqzero}
\Phi_g(z_0)=0,
\end{equation}
so the nonnegative function $z\mapsto\Phi_g(z)$ attains its minimum value $0$ at $z_0$. This function being also superharmonic on $\D$, it obeys the minimum principle (equivalently, the subharmonic function $-\Phi_g$ obeys the maximum principle), and this means that \eqref{eqn:phigeqzero} is only possible if $\Phi_g(z)=0$ for all $z\in\D$. Consequently, we have 
$g=\const$ as desired.\qed

\bigskip\noindent{\it Proof of Theorem \ref{thm:extpoilip}.} This is similar to the preceding proof, except for some routine changes. Essentially, these amount to dividing every term in \eqref{eqn:pmphiphi} by $(1-|z|)^{2\al}$ and proceeding as before, always using the norm $\|\cdot\|_{G,\al}$ in place of $\|\cdot\|_G$. We omit further details.\qed

\bigskip\noindent{\it Proof of Theorem \ref{thm:nonextpoibmo}.} (a) Suppose that $\ph\in L^\infty$ is a non-unimodular function with $\|\ph\|_\infty=1$. (In other words, $\ph$ is a unit-norm $L^\infty$ function with nonconstant modulus.) We know from Theorem \ref {thm:imprax} that there exist a Blaschke product $B$ and a $G$-extremal function $f\in L^\infty$ such that 
\begin{equation}\label{eqn:popka}
B\ph=f.
\end{equation}
Thus, in particular, 
\begin{equation}\label{eqn:ginftyone}
\|f\|_G=\|f\|_\infty=1.
\end{equation}
Furthermore, letting 
\begin{equation}\label{eqn:defggg}
g:=1-|\ph|,
\end{equation}
we define 
\begin{equation}\label{eqn:defphionetwo}
\ph_1:=\ph+g\quad\text{\rm and}\quad\ph_2:=\ph-g,
\end{equation}
so that $\ph_1,\ph_2\in L^\infty$ and
\begin{equation}\label{eqn:zadok}
\ph=\f12\left(\ph_1+\ph_2\right).
\end{equation}
Together with \eqref{eqn:popka} this yields
\begin{equation}\label{eqn:fbphibphi}
f=\f12\left(B\ph_1+B\ph_2\right).
\end{equation}
Now, one quickly checks for $j=1,2$ that $\|\ph_j\|_\infty\le1$, and so 
\begin{equation}\label{eqn:peredok}
\|B\ph_j\|_G\le\|B\ph_j\|_\infty=\|\ph_j\|_\infty\le1.
\end{equation}
Consequently, 
\begin{equation}\label{eqn:zhopka}
B\ph_j\in{\rm ball}\,({\bmo}),\qquad j=1,2.
\end{equation}
Finally, we observe that the functions $B\ph_1$ and $B\ph_2$ represent two distinct elements of $\bmo$, meaning that 
\begin{equation}\label{eqn:difnoncon}
B\ph_1-B\ph_2\ne\const.
\end{equation}
To see why, it suffices to note that
\begin{equation}\label{eqn:stolik}
|B\ph_1-B\ph_2|=|\ph_1-\ph_2|=2|g|=2(1-|\ph|)
\end{equation}
and recall that $|\ph|$ is nonconstant. This shows that 
\begin{equation}\label{eqn:stulik}
|B\ph_1-B\ph_2|\ne\const,
\end{equation}
whence \eqref{eqn:difnoncon} is immediate. In light of this last observation, we readily deduce from \eqref{eqn:fbphibphi} and \eqref{eqn:zhopka} that $f$ is a non-extreme point of ${\rm ball}\,({\bmo})$. 

\medskip (b) To prove this part, which is an \lq\lq analytic version" of (a), we only need to make a couple of adjustments in the above reasoning. Suppose that $\ph\in H^\infty$ is a function with $\|\ph\|_\infty=1$ that obeys \eqref{eqn:logintegconv}. Using Theorem \ref{thm:imprax} as before (this time, in conjunction with the remark that follows it), we arrive at \eqref{eqn:popka} with a suitable Blaschke product $B$ and a $G$-extremal function $f\in H^\infty$. Thus, \eqref{eqn:ginftyone} is again fulfilled. 
\par Next, we write $g$ for the outer function with modulus $1-|\ph|$ on $\T$ (condition \eqref{eqn:logintegconv} makes this meaningful) and then define $\ph_1$ and $\ph_2$ as in \eqref{eqn:defphionetwo}, but with the current $g$ in place of its namesake given by \eqref{eqn:defggg}. The functions $\ph_j$ ($j=1,2$) built in this way are in $H^\infty$, not just in $L^\infty$, and the formulas \eqref{eqn:zadok}--\eqref{eqn:stulik} remain valid, as does the text that accompanies them. Moreover, \eqref{eqn:zhopka} can now be written in the form 
$$B\ph_j\in{\rm ball}\,({\bmoa}),\qquad j=1,2.$$
Combining this with \eqref{eqn:fbphibphi}, while taking \eqref{eqn:difnoncon} into account, we see that $f$ is a non-extreme point of ${\rm ball}\,({\bmoa})$. 

\medskip (c) The conclusion that $B\ph$ is a non-extreme point of the unit ball in $\bmo$, or in $\bmoa$, is immediate if $\|B\ph\|_G<1$. The remaining cases have already been settled when proving parts (a) and (b) above. Indeed, no specific assumptions on the Blaschke product $B$ (other than the condition $\|B\ph\|_G=1$) were imposed there. Thus, the desired conclusion actually holds for an arbitrary $B$, and we are done.\qed

\bigskip\noindent{\it Proof of Theorem \ref{thm:charint}.} (i)$\implies$(ii). If $h$ is inner and $B$ is a nonconstant Blaschke product, then $Bh$ is an inner function that has a zero in $\D$. Therefore, by Corollary \ref{cor:ginn}, this last function is both norm-attaining and $G$-extremal; that is, $Bh\in\bmona$ and 
$$\|Bh\|_G=\|Bh\|_\infty=1.$$
The fact that $Bh$ is an extreme point of ${\rm ball}\,({\bmo})$ is now ensured by Theorem \ref{thm:extpoibmo}. 

\smallskip(ii)$\implies$(iii). This is obvious. 

\smallskip(iii)$\implies$(i). Assuming that $Bh$ is an extreme point of ${\rm ball}\,({\bmo})$ for some Blaschke product $B$, we readily deduce that $h$ is unimodular (and hence inner) by applying part (c) of Theorem \ref{thm:nonextpoibmo}, with $h$ in place of $\ph$. The proof is complete.\qed

\section{Example: a $G$-extremal outer function}

Here we give an explicit example of a nonnegative function $\eta\in L^\infty$, with 
\begin{equation}\label{eqn:suslik}
\int_\T\log\eta\,dm>-\infty,
\end{equation}
for which the corresponding outer function 
$$\mathcal O_\eta(z):=\exp\left\{\int_\T\f{\ze+z}{\ze-z}\log\eta(\ze)\,dm(\ze)\right\},\qquad z\in\D,$$
is $G$-extremal. Our $\eta$ will have $\|\eta\|_\infty=1$, and we shall exhibit a sequence $\{z_k\}\subset\D$ with the properties that
\begin{equation}\label{eqn:vilochka}
\mathcal P\eta(z_k)\to1
\end{equation}
and
\begin{equation}\label{eqn:nozhichek}
\mathcal P\left(\log\f1\eta\right)(z_k)\to\infty.
\end{equation}
Once this is achieved, Theorem \ref{thm:charoutergextr} will guarantee that $\mathcal O_\eta$ is indeed $G$-extremal. 
\par Now let $\{l_k\}$, $\{\eps_k\}$ and $\{d_k\}$ be three sequences of real numbers from the interval $(0,1)$, all of them tending to $0$ and such that 
$$\sum_{k=1}^\infty l_k\le\pi,\qquad\sum_{k=1}^\infty l_k\log\f1{\eps_k}<\infty,\qquad\f{d_k}{l_k}\to0,
\qquad\f{d_k}{l_k}\log\f1{\eps_k}\to\infty.$$
For instance, we may take 
$$l_k=\eps_k=2^{-k}\quad\text{\rm and}\quad d_k=k^{-1/2}2^{-k}.$$
We then define the numbers 
$$0=\al_1<\be_1<\al_2<\be_2<\dots$$
inductively by the relations
$$\al_{k+1}-\be_k=\be_k-\al_k=l_k\qquad(k=1,2,\dots)$$
and consider the (pairwise disjoint) interlacing arcs
$$I_k:=\left\{e^{it}:\,\al_k\le t<\be_k\right\}\quad\text{\rm and}\quad J_k:=\left\{e^{it}:\,\be_k\le t<\al_{k+1}\right\}.$$
Finally, we put $\mathcal J:=\bigcup_{k=1}^\infty J_k$ and introduce the function
$$\eta:=\chi_{\T\setminus\mathcal J}+\sum_{k=1}^\infty\eps_k\chi_{J_k},$$
which is obviously positive and satisfies $\|\eta\|_\infty=1$. 

\par To check that this $\eta$ actually does the job, we first verify \eqref{eqn:suslik}. In fact, we have 
$$\int_\T\log\f1\eta\,dm=\sum_{k=1}^\infty\int_{J_k}\log\f1\eta\,dm=\f1{2\pi}\sum_{k=1}^\infty l_k\log\f1{\eps_k}<\infty.$$
Next, we define the sequence $\{z_k\}\subset\D$ by 
$$z_k=(1-d_k)e^{it_k}\qquad(k=1,2,\dots),$$
where $t_k:=(\al_k+\be_k)/2$, and we are going to show that this makes \eqref{eqn:vilochka} and \eqref{eqn:nozhichek} true. Indeed, since $\eta|_{I_k}=1$ for each $k$, it follows that 
\begin{equation}\label{eqn:toporik}
\mathcal P\eta(z_k)\ge\int_{I_k}\eta\,d\om_{z_k}=\int_{I_k}d\om_{z_k}=\om_{z_k}(I_k).
\end{equation}
Condition \eqref{eqn:vilochka} is now implied by \eqref{eqn:toporik}, coupled with the facts that $\mathcal P\eta(z_k)\le1$ and 
\begin{equation}\label{eqn:tyapka}
\om_{z_k}(I_k)\to1.
\end{equation} 
(To verify \eqref{eqn:tyapka}, note that the radial projection of $z_k$ is the midpoint of $I_k$ and recall that $d_k/l_k\to0$.) Furthermore, we have $\eta|_{J_k}=\eps_k$ for each $k$, and so
\begin{equation}\label{eqn:motyga}
\mathcal P\left(\log\f1\eta\right)(z_k)\ge\int_{J_k}\log\f1\eta\,d\om_{z_k}
=\om_{z_k}(J_k)\cdot\log\f1{\eps_k}.
\end{equation}
A calculation shows that $\om_{z_k}(J_k)$ is comparable to $d_k/l_k$, and since 
$$\f{d_k}{l_k}\log\f1{\eps_k}\to\infty$$
by assumption, we readily deduce \eqref{eqn:nozhichek} from \eqref{eqn:motyga}. The required properties are thereby established. 

\par We finally mention, in connection with the above example, that our current construction is somewhat reminiscent of that in \cite[Section 3]{DRMI}, where a similar method was used to produce a non-$\bmoa$ outer function with $\bmo$ modulus.

\end{document}